\documentclass[10pt,twocolumn,abstract,BCOR=0pt,DIV=calc]{scrartcl}

\usepackage{amsmath}
\usepackage{amssymb}
\usepackage{amsthm}
\usepackage[colorlinks = true, linkcolor = blue, citecolor = blue]{hyperref}
\usepackage[nameinlink]{cleveref}
\usepackage{graphicx}

\setcapindent{0pt}
\setkomafont{captionlabel}{\bfseries}

\newtheorem{assumption}{Assumption}
\newtheorem{lemma}{Lemma}
\newtheorem{proposition}{Proposition}
\newtheorem{remark}{Remark}
\newtheorem{theorem}{Theorem}

\crefname{assumption}{Assumption}{Assumptions}
\crefname{remark}{Remark}{Remarks}

\crefname{equation}{}{}
\crefname{figure}{Figure}{Figures}

\crefrangelabelformat{subequation}{(#3#1#4--#5\crefstripprefix{#1}{#2}#6)}

\begin{document}

\title{Nonlocal Nonholonomic Source Seeking \\ Despite Local Extrema}

\author{Raik Suttner,
\thanks{This research was supported by the German Research Foundation DFG, project number DA~767/13-1.}
\thanks{Raik Suttner is with the Institute of Mathematics, University of Wuerzburg, Wuerzburg, Germany (e-mail: raik.suttner@mathe\-matik.uni-wuerzburg.de)}
$\quad$ Miroslav Krsti\'{c}
\thanks{Miroslav Krsti\'{c} is with the Department of Mechanical and Aero\-space Engineering, University of California, San Diego, USA (e-mail: krstic@ucsd.edu)}}

\maketitle


\begin{abstract}
In this paper, we investigate the problem of source seeking with a unicycle in the presence of local extrema. Our study is motivated by the fact that most of the existing source seeking methods follow the gradient direction of the signal function and thus only lead to local convergence into a neighborhood of the nearest local extremum. So far, only a few studies present ideas on how to overcome local extrema in order to reach a global extremum. None of them apply to second-order (force- and torque-actuated) nonholonomic vehicles. We consider what is possibly the simplest conceivable algorithm for such vehicles, which employs a constant torque and a translational/surge force in proportion to an approximately differentiated measured signal. We show that the algorithm steers the unicycle through local extrema towards a global extremum. In contrast to the previous extremum-seeking studies, in our analysis we do not approximate the gradient of the objective function but of the objective function's local spatial average. Such a spatially averaged objective function is expected to have fewer critical points than the original objective function. Under suitable assumptions on the averaged objective function and on sufficiently strong translational damping, we show that the control law achieves practical uniform asymptotic stability and robustness to sufficiently weak measurement noise and disturbances to the force and torque inputs.
\end{abstract}


\section{Introduction}\label{Section:1}
The problem of \emph{source seeking} arises when an autonomous agent has to locate a point where an unknown scalar signal attains an extreme value. For instance, the autonomous agent could be a wheeled robot or a drone and the unknown scalar signal could be given by the strength of an electromagnetic field or the concentration of a chemical substance. It is assumed that the agent is equipped with a suitable sensor so that it can measure the signal at its current position. The spatial distribution of the signal is unknown. Moreover, the agent can neither determine its own position or orientation with respect to a reference frame nor can it store data to compare measurement results at different time instances. Only real-time measurements of the signal values can be used for a feedback law. For many real-world signals, it is reasonable to expect that the signal attains an unknown global extreme value at some unknown point, called the \emph{source}. In this sense, source seeking can be seen as a global optimization problem. Here, we follow the usual convention in optimization and focus (without loss of generality) on minimization. 

Over the past two decades, many different source-seeking methods for various agent models have been proposed. This includes methods for single-integrator points \cite{Duerr20112,Stankovic2009}, velocity-actuated unicycles in the plane \cite{Liu20102,Duerr2017}, and velocity-actuated nonholonomic models in the three-dimensional space \cite{Cochran20092,Matveev2014}. However, velocity actuation considers only vehicle kinematics. Wheeled robots and drones are mechanical systems with inertia, and, therefore, second-order dynamic models, actuated by forces and torques, are more realistic descriptions than first-order kinematic models actuated by velocities. There exist source-seeking methods for double-integrator points \cite{Zhang20072,Michalowsky2014}, force- and torque-controlled slider models \cite{Scheinker20183,Suttner20221}, and force- and torque-controlled unicycles \cite{Suttner20192}. In the present paper, we focus on the latter model. The unicycle has two inputs: A force to control the translational motion and a torque to control the rotational motion.

As indicated above, the objective of source seeking is, ideally, global optimization. That is, the goal is to find the point where the unknown signal attains its global minimum value; under the tacit assumption that such a point exists and is unique. We also assume that the signal is given by a sufficiently smooth real-valued function, which then plays the role of an objective function. Many of the existing source-seeking methods steer the agent approximately into the negative gradient direction of the objective function. This approach works very well as long as the objective function has no other critical points besides the global minimum. Otherwise, it can happen that the agent gets stuck in an undesired configuration. In general, a gradient-based algorithm only leads to convergence towards the set of critical points but not to the global minimum. The same problem also occurs in the context of \emph{extremum-seeking control} for general non-linear systems. So far, relatively little is known about systematic procedures to overcome undesired local extrema. The intention of the present paper is to offer a partial solution, at least for the particular problem of source seeking with a dynamic unicycle.

\emph{Survey of literature.} Global extremum seeking in the presence of local extrema is certainly a very difficult problem. Many papers, such as \cite{Nesic2006,Wang2016,Bhattacharjee2021}, address this topic only in terms of qualitative discussions and simulation studies. Quantitative statements about the ability of a method to overcome local extrema are rather difficult to derive. There are, however, some attempts to provide sufficient conditions under which a method can pass through a local extremum; for instance in \cite{Tan2006,Tan2009} for a steady-state input-output map of one variable. The proposed method in \cite{Tan2006,Tan2009} employs a periodic perturbation signal with tunable amplitude. Under suitable assumptions, it is shown in \cite{Tan2006,Tan2009} that if the amplitude is sufficiently large, then the input can pass through local extrema of the input-output function. In the present paper, we will interpret this behavior in terms of the gradient of an averaged input-output function. Such an averaged input-output function is expected to have fewer local extrema than the original input-output function. The control law from \cite{Tan2006,Tan2009} is extended in \cite{Ye2020} to solve a non-convex social cost optimization problem by a network of agents. In this context, we also mention the perturbation-based extremum-seeking methods for open-loop unstable kinematic systems in \cite{Wildhagen2018,Suttner20222}. The investigations in \cite{Wildhagen2018,Suttner20222} lead to a similar interpretation as for the method in \cite{Tan2006,Tan2009}. It is conjectured in \cite{Wildhagen2018,Suttner20222} that the employed periodic perturbation signals give approximate access to the gradient of an averaged objective function. However the investigations in \cite{Wildhagen2018,Suttner20222} do not reveal explicitly how the suspected averaged objective function is related to the original objective function. A completely different approach to global optimization is taken in \cite{EsmaeilZadehAzar2010,EsmaeilZadehAzar2011} in the form of multi-unit extremum-seeking for one- and two-dimensional single-integrator points. However, an implementation of the method in \cite{EsmaeilZadehAzar2010,EsmaeilZadehAzar2011} requires that the objective function  be measured at different positions at the same time. The methods in \cite{Nesic2013,Khong2013} take advantage of the fact that some established discrete-time derivative-free algorithms, such as the Piyavskii–Shubert algorithm \cite{Piyavskii1972,Shubert1972}, have the ability to find the global minimum of a function (within a prescribed compact set). Using a zero-order hold, one can apply these discrete-time algorithms to a continuous-time system for the purpose of global extremum-seeking (under the assumption of a suitable steady-state input-output behavior of the underlying control system). Finally, we note that none of the above methods is suitable for source-seeking with a dynamic unicycle in the presence of local extrema.

\emph{Underlying idea.} The results in the above-mentioned papers \cite{Tan2006,Tan2009,Wildhagen2018,Suttner20222} suggest that oscillations of the sensor with sufficiently large amplitudes could be beneficial for our goal of source-seeking in the presence of local extrema, because these oscillations might give access to an averaged signal function with fewer undesired critical points than the original signal function. To generate a suitable oscillatory motion, we demand that the sensor is attached to the unicycle at a point with distance $r$ from the center of the vehicle; the position of the vehicle center is denoted by $p$. Such a displaced sensor was first proposed in \cite{Cochran2009}, and then also used in \cite{Ghods2010,Raisch2017}, for the kinematic unicycle model, but these papers do not address the problem of source-seeking in the presence of local extrema. If the sensor is mounted in the described way, then a rotational motion of the unicycle results in a circular motion of the sensor. Note that a circular motion of the sensor requires that the translational velocity is slow compared to rotational velocity. For this reason, we demand sufficiently strong linear damping of the translational velocity. A suitable design of our feedback law leads to an integration of the measured values of the signal function $\psi$. After one full rotation of the unicycle, we get the integral of $\psi$ over the circle with radius $r$ centered at $p$. In this situation, we can apply the \emph{divergence theorem} (a consequence of \emph{Stokes' theorem}) to write the integral over the circle in terms of an integral over the disk. The latter integral leads us to the definition of the averaged signal function $\Psi_r$ at $p$, namely as the integral of $\psi$ over the disk of radius $r$ centered at $p$. Such an application of the divergence theorem is also used in some studies on nonsmooth optimization to obtain a smoothed objective function; see, e.g., \cite{Gupal1977,Mayne1984,Flaxman2005}. If the sensor distance $r$ is sufficiently large, then one can expect that the averaged signal function $\Psi_r$ has fewer undesired critical points than the original signal function $\psi$. Moreover, a global minimum of $\Psi_r$ is expected to be close to the global minimum of $\psi$. Consequently, the function $\Psi_r$ is better suited for global minimization by gradient descent than $\psi$. We will show that the unicycle under the proposed control law follows approximately the negative gradient direction of $\Psi_r$. As a result, we can expect that the vehicle center converges into a small neighborhood of the source.

\emph{Novelties/Contributions.} Since our method only requires a uniform rotation of the unicycle with moderate rotational velocity, it is rather easy to implement and also less energy consuming than many other source-seeking methods with large-amplitude high-frequency inputs, such as in \cite{Zhang20071,Cochran2009,Duerr20112}, which generate an oscillatory motion of the unicycle. However, this goes along with the disadvantage that the speed of convergence of our method is rather slow. In contrast to many other studies on source-seeking, our theoretical analysis takes $L_\infty$-input disturbances and $L_\infty$-measurement errors into account. This proves that our approach grantees a certain degree of robustness with respect to sufficiently weak disturbances. As indicated in the previous paragraph, we can provide a simple formula for the averaged signal function. Our main stability result only involves assumptions on the averaged signal function. Those assumptions can be satisfied even if the original signal function has undesired local minima. Under these mild assumptions, we prove practical uniform asymptotic stability for the closed-loop system in the presence of sufficiently weak disturbances. The proof of our stability result involves a suitable first-order averaging analysis for the closed-loop system. The averaging analysis reveals that the position of the unicycle is driven approximately into the negative gradient direction of the averaged signal function if a certain control parameter is sufficiently small. The quality of approximation improves with shrinking parameter. This in turn leads to the effect that robustness and stability properties of the gradient system carry over to the approximating closed-loop system for sufficiently small parameter values. Convergence to the global minimum in the presence of local minima is illustrated by means of a numerical example.


\section{Problem statement and control law}\label{Section:2}
In the first step, we introduce the dynamic equations for a second-order unicycle. The unicycle is a mechanical system with body mass $m>0$ and moment of inertia $J>0$. The configuration of the unicycle at time $t\in\mathbb{R}$ is determined by its current position $p(t)\in\mathbb{R}^2$ and its current orientation $o(t)\in\mathbb{S}$, where $\mathbb{S}$ is the unit circle centered at the origin. We use the variables $v(t)\in\mathbb{R}$ and $\omega(t)\in\mathbb{R}$ to describe the translational and rotational velocity, respectively. Let $k$ and $\kappa$ be arbitrary positive real constants to describe the impact of linear velocity-dependent damping. Finally, the letter $R$ denotes the $(2\times{2})$-rotation matrix for counter-clockwise rotation by angle $\pi/2$. Then $Ro(t)$ is a tangent vector to $\mathbb{S}$ at $o(t)$. We consider the second-order unicycle model%
\begin{subequations}\label{eq:01}%
\begin{align}
\dot{p}(t) & \ = \ v(t)\,o(t), \label{eq:01:a} \allowdisplaybreaks \\
\dot{o}(t) & \ = \ \omega(t)\,Ro(t), \label{eq:01:b} \allowdisplaybreaks \\
m\,\dot{v}(t) & \ = \ - \tfrac{k}{\varepsilon}\,v(t) + \mathrm{u}_{\text{t}}(t) + d_{\text{t}}(t), \label{eq:01:c} \allowdisplaybreaks \\
J\,\dot{\omega}(t) & \ = \ - \kappa\,\omega(t) + \mathrm{u}_{\text{r}}(t) + d_{\text{r}}(t), \label{eq:01:d}
\end{align}%
\end{subequations}%
where $\mathrm{u}_{\text{t}}$, $\mathrm{u}_{\text{r}}$ are real-valued inputs and $d_{\text{t}},d_{\text{r}}\in{L_\infty}$ are input disturbances. Here, $L_\infty$ denotes the space of essentially bounded and measurable real-valued functions on the real line. The essential supremum norm on $L_\infty$ is denoted by $\|\cdot\|$. In some texts, e.g. in \cite{BulloBook}, the dynamic unicycle model \cref{eq:01} is called \emph{rolling disk model}.

In the vehicle model~\cref{eq:01}, we assume the existence of linear velocity-dependent damping. The presence of dissipation prevents an unbounded increase of the total energy (source seeking control for fully actuated systems in the absence of dissipation is studied in \cite{Suttner20232}). The constants $k/\varepsilon$ and $\kappa$ in \cref{eq:01:c,eq:01:d} represent the translational and rotational damping coefficients, respectively, where $\varepsilon$ is a parameter that is small relative to $k,\kappa$. Thus, we assume that the vehicle has been designed to exhibit higher damping in translation than in rotation, i.e., to promote rotation over translation. This is crucial in our control design and analysis.

\begin{figure}%
\centering\includegraphics{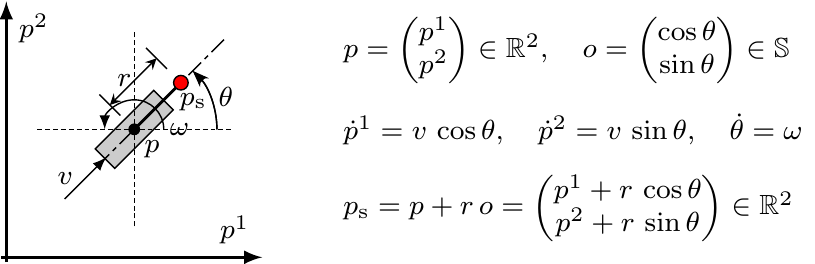}%
\caption{Notation used in the model of vehicle sensor and center dynamics. The position of the sensor is denoted by $p_{\text{s}}$.}%
\label{Figure:1}%
\end{figure}

We assume that the unicycle is equipped with a suitable sensor so that it can measure the value of an unknown scalar signal. The signal is assumed to be given by a smooth real-valued function $\psi$ on $\mathbb{R}^2$. We do \emph{not} assume that $\psi$ is analytically known. In particular, the gradient of $\psi$ is an unknown quantity. The proposed method requires that the sensor is non-collocated with the center of the vehicle. To make it precise, we assume that the sensor is mounted at the point $p(t)+ro(t)$ with fixed distance $r>0$ from the vehicle center $p(t)$ into the direction $o(t)$ of the current alignment. Then, a rotational motion of the unicycle results in a circular motion of the sensor. This way, the sensor can explore the unknown signal in a neighborhood of the current position of the vehicle center. A suitable design of our feedback law will allow us to extract gradient information from the circular motion of the sensor. A measurement of $\psi$ at the current sensor position $p(t)+ro(t)\in\mathbb{R}^2$ at time $t\in\mathbb{R}$ results in the real value%
\begin{equation}\label{eq:02}
\hat{\mathrm{y}}(t) \ = \ \psi(p(t)+r\,o(t)) + d_{\text{s}}(t),
\end{equation}%
where the function $d_{\text{s}}\in{L_\infty}$ describes disturbances of the signal measurements. One may view \cref{eq:02} as the noise-corrupted output of the unicycle system. We are interested in an output feedback law that asymptotically stabilizes the unicycle at positions where $\psi$ attains a minimum value. The diagram in \Cref{Figure:1} depicts the position, heading, angular and forward velocities for the center and sensor.%

We propose the control law%
\begin{subequations}\label{eq:03}%
\begin{align}
\mathrm{u}_{\text{t}}(t) & \ = \ -\lambda\,\big(\hat{\mathrm{y}}(t)-\eta(t)\big), \label{eq:03:a} \allowdisplaybreaks \\
\mathrm{u}_{\text{r}}(t) & \ = \ \tau_\ast, \label{eq:03:b}
\end{align}%
\end{subequations}%
where $\lambda$, $\tau_\ast$ are arbitrary positive constants and the high-pass filter%
\begin{equation}\label{eq:04}
\dot{\eta}(t) \ = \ \varepsilon\,h\,\big(\hat{\mathrm{y}}(t)-\eta(t)\big)
\end{equation}%
with filter input $\hat{\mathrm{y}}$ and filter output $\hat{\mathrm{y}}-\eta$ removes the DC bias in \cref{eq:02}, where $h$ is a positive constant. The inclusion of $\varepsilon$ into the filter gain $\varepsilon\,h$, which is inversely proportional to the translational damping coefficient $k/\varepsilon$, does not mean that the value of the translational damping is known but only that the damping's order of magnitude is known. Setting the filter's pole in inverse proportion to  translational damping slows the filter relative to the rate of rotation induced by the constant torque $\tau_\ast$ and low rotational damping $-\kappa\,\omega(t)$ in \cref{eq:01:d}. The rotational velocity $\omega(t)$ converges into a small neighborhood of%
\begin{equation}\label{eq:05}
\omega_\ast \ := \ \tau_\ast/\kappa \ > \ 0
\end{equation}%
if the magnitude $\|d_{\text{r}}\|$ of the input disturbance $d_{\text{r}}\in{L_\infty}$ in~\cref{eq:01:d} is sufficiently small.

For what follows, we summarize the closed-loop system by substituting \cref{eq:03}, \cref{eq:04} into \cref{eq:01}, \cref{eq:02} to get the system%
\begin{subequations}\label{eq:06}%
\begin{align}
\dot{p}(t) & \ = \ v(t)\,o(t), \label{eq:06:a} \allowdisplaybreaks \\
\dot{o}(t) & \ = \ \omega(t)\,Ro(t), \label{eq:06:b} \allowdisplaybreaks \\
m\,\dot{v}(t) & \ = \ -\lambda\big(\psi(p(t)+r\,o(t)) - \eta(t)\big) - \tfrac{1}{\varepsilon}\,k\,v(t) \label{eq:06:c} \allowdisplaybreaks \\
& \qquad - \lambda\,d_{\text{s}}(t) + d_{\text{t}}(t), \label{eq:06:d} \allowdisplaybreaks \\
J\,\dot{\omega}(t) & \ = \ - \kappa\,(\omega(t) - \omega_\ast) + d_{\text{r}}(t), \label{eq:06:e} \allowdisplaybreaks \\
\dot{\eta}(t) & \ = \ -\varepsilon\,h\,\eta(t) + \varepsilon\,h\,\psi(p(t)+r\,o(t)) + \varepsilon\,h\,d_{\text{s}}(t) \label{eq:06:f}
\end{align}%
\end{subequations}%
on the state manifold $M:=\mathbb{R}^2\times\mathbb{S}\times\mathbb{R}\times\mathbb{R}\times\mathbb{R}$.


\section{Averaging}\label{Section:3}
In this section, we will see that the closed-loop system \cref{eq:06} approximates the behavior of a certain averaged system if the positive control parameters $\varepsilon$ is sufficiently small. For the sake of brevity and clarity, several technical details are omitted here. The reader is referred to the \hyperlink{appendixChapter}{Appendix} for a more precise averaging analysis.


\subsection{Averaging of the rotational motion}\label{Section:3.1}
To reveal the impact of the rotational motion with rotational velocity $\omega_\ast$, we carry out the change of variables%
\begin{subequations}\label{eq:07}%
\begin{align}
\tilde{p}(\tau) & \ = \ p(\tfrac{\tau}{\varepsilon}), \label{eq:07:a} \allowdisplaybreaks \\
\tilde{o}(\tau) & \ = \ \exp(-\omega_\ast\,\tfrac{\tau}{\varepsilon}\,R)\,o(\tfrac{\tau}{\varepsilon}), \allowdisplaybreaks \\
\tilde{v}(\tau) & \ = \ v(\tfrac{\tau}{\varepsilon}) + \varepsilon\,\tfrac{\lambda}{k}\,\big(\psi(p(\tfrac{\tau}{\varepsilon})+r\,o(\tfrac{\tau}{\varepsilon}))-\eta(\tfrac{\tau}{\varepsilon})\big), \allowdisplaybreaks \\
\tilde{\omega}(\tau) & \ = \ \omega(\tfrac{\tau}{\varepsilon}) - \omega_\ast, \allowdisplaybreaks \\
\tilde{\eta}(\tau) & \ = \ \eta(\tfrac{\tau}{\varepsilon}),
\end{align}%
\end{subequations}%
where $\exp\colon\mathbb{R}^{2\times{2}}\to\mathbb{R}^{2\times{2}}$ is the matrix exponential function. If the control parameter $\varepsilon>0$ is small, then, in the new time scale $\tau$, the unicycle is rotating fast while its position is varying slowly. This means that the sensor performs a fast circular motion around the current position of the unicycle. A direct computation shows that the closed-loop system \cref{eq:06} in the variables \cref{eq:07} can be equivalently written as%
\begin{subequations}\label{eq:08}%
\begin{align}
\dot{\tilde{p}}(\tau) & \ = \ -\tfrac{\lambda\,r}{2\,k}\,\tilde{G}_r(\tfrac{\tau}{\varepsilon},\tilde{p}(\tau),\tilde{o}(\tau),\tilde{\eta}(\tau)) \label{eq:08:a} \allowdisplaybreaks \\
& \qquad + \tfrac{1}{\varepsilon}\,\tilde{v}(\tau)\,\exp(\omega_\ast\,\tfrac{\tau}{\varepsilon}\,R)\tilde{o}(\tau) \label{eq:08:b} \allowdisplaybreaks \\
\dot{\tilde{o}}(\tau) & \ = \ \tfrac{1}{\varepsilon}\,\tilde{\omega}(\tau)\,R\tilde{o}(\tau), \label{eq:08:c} \allowdisplaybreaks \\
m\,\dot{\tilde{v}}(\tau) & \ = \ -\tfrac{1}{\varepsilon^2}\,k\,\tilde{v}(\tau) + \tilde{Q}(\tfrac{\tau}{\varepsilon},\tilde{x}(\tau)) + \tfrac{1}{\varepsilon}\,d_{\text{t}}(\tfrac{\tau}{\varepsilon}) \label{eq:08:d} \allowdisplaybreaks \\
& \qquad - \lambda\,\big(\tfrac{1}{\varepsilon} + \varepsilon\,\tfrac{h\,m}{k}\big)\,d_{\text{s}}(\tfrac{\tau}{\varepsilon}) + \varepsilon\,\tilde{R}(\tfrac{\tau}{\varepsilon},\tilde{x}(\tau)), \label{eq:08:e} \allowdisplaybreaks \\
J\,\dot{\tilde{\omega}}(\tau) & \ = \ -\tfrac{1}{\varepsilon}\,\kappa\,\tilde{\omega}(\tau) + \tfrac{1}{\varepsilon}\,d_{\text{r}}(\tfrac{\tau}{\varepsilon}), \label{eq:08:f} \allowdisplaybreaks \\
\dot{\tilde{\eta}}(\tau) & \ = \ -h\,\tilde{\eta}(\tau) + h\,d_{\text{s}}(\tfrac{\tau}{\varepsilon}) \label{eq:08:g} \allowdisplaybreaks \\
& \qquad + h\,\psi\big(\tilde{p}(\tau)+r\,\exp(\omega_\ast\,\tfrac{\tau}{\varepsilon}\,R)\tilde{o}(\tau)\big), \label{eq:08:h}
\end{align}%
\end{subequations}%
where we write $\tilde{x}=(\tilde{p},\tilde{o},\tilde{v},\tilde{\omega},\tilde{\eta})$, the $\mathbb{R}^2$-valued map $\tilde{G}_r$ on $\mathbb{R}\times\mathbb{R}^2\times\mathbb{S}\times\mathbb{R}$ is defined by%
\begin{subequations}\label{eq:09}%
\begin{align}
\tilde{G}_r&(t,p,o,\eta) \ := \ \allowdisplaybreaks \\
& \tfrac{2}{r}\big(\psi(p+r\,\exp(\omega_\ast{t}\,R)o)-\eta\big)\,\exp(\omega_\ast{t}\,R)o,
\end{align}%
\end{subequations}%
and $\tilde{Q}(t,\cdot)$, $\tilde{R}(t,\cdot)$ are certain smooth real-valued functions on the state manifold $M$ which are $(2\pi/\omega_\ast)$-periodic in $t$. If the disturbances $d_{\text{r}}$, $d_{\text{s}}$, and $d_{\text{t}}$ are sufficiently weak, then the transformed velocities $\tilde{v} $ and $\tilde{\omega}$ in \cref{eq:08:d,eq:08:e,eq:08:f} converge fast into small neighborhoods of the origin. After short transients of $\tilde{v} $ and $\tilde{\omega}$, the position $\tilde{p}$, the transformed orientation $\tilde{o}$, and the high-pass filter $\tilde{\eta}$ in \cref{eq:08:a,eq:08:b,eq:08:c,eq:08:g,eq:08:h} can be shown to be moderately time-varying compared to the strongly varying periodic rotation $\tau\mapsto\exp(\omega_\ast\,\frac{\tau}{\varepsilon}\,R)$ with rotational velocity $\omega_\ast/\varepsilon$. In this situation, it is reasonable to expect that the impact of the high-frequency vector field on the right-hand side of \cref{eq:08:a} is approximately the same as its time average. One can easily check that, for all fixed $\tilde{p}\in\mathbb{R}^2$, $\tilde{o}\in\mathbb{S}$, and $\tilde{\eta}\in\mathbb{R}$, the average of $t\mapsto\tilde{G}_r(t,\tilde{p},\tilde{o},\tilde{\eta})$ is given by%
\begin{equation}\label{eq:10}
\tfrac{\omega_\ast}{2\pi}\int_{0\qquad}^{\frac{2\pi}{\omega_\ast}}\!\!\!\!\!\!\tilde{G}_r(t,\tilde{p},\tilde{o},\tilde{\eta})\,\mathrm{d}t \ = \ \bar{G}_r(\tilde{p}),
\end{equation}%
where the vector field $\bar{G}_r$ on $\mathbb{R}^2$ is defined by%
\begin{equation}\label{eq:11}
\bar{G}_r(\bar{p}) \ := \ \tfrac{1}{\pi}\int_{\mathbb{S}}\tfrac{1}{r}\,\psi(\bar{p}+r\,o)\,o\,\mathrm{d}\lambda_\mathbb{S}(o)
\end{equation}%
and $\lambda_{\mathbb{S}}$ denotes the standard volume measure on $\mathbb{S}$. Recall that \cref{eq:07:a} is just a change of the time scale. Thus, our vague averaging argument indicates that, for sufficiently weak disturbances $d_{\text{r}}$, $d_{\text{s}}$, and $d_{\text{t}}$ and a sufficiently small control parameter $\varepsilon>0$, the solutions of \cref{eq:06:a} approximate the solutions of a system of the form%
\begin{equation}\label{eq:12}
\dot{\bar{p}}(t) \ = \ -\varepsilon\,\tfrac{\lambda\,r}{2\,k}\,\bar{G}_r(\bar{p}(t)) + \varepsilon\,\bar{d}(t),
\end{equation}%
where the ``disturbance'' $t\mapsto\bar{d}(t)\in\mathbb{R}^2$ originates from the contribution of \cref{eq:08:b} in the original time scale $t$. A precise approximation result is given by \Cref{Proposition:1} in the \hyperlink{appendixChapter}{Appendix}. In the next subsection, we will see that the averaged equation \cref{eq:12} is a disturbed negative gradient system of an averaged signal function.


\subsection{Averaged signal function}\label{Section:3.2}
The vector field $\bar{G}_r$ in \cref{eq:11} is crucial for our objective of source seeking in the presence of undesired local minima. Let $P$ be an arbitrary non-empty open subset of $\mathbb{R}^2$; for instance, $P=\mathbb{R}^2$. For every positive real number $r$, define a real-valued function $\Psi_r$ on $P$ by%
\begin{equation}\label{eq:13}
\Psi_r(\bar{p}) \ := \ \tfrac{1}{\pi}\int_{\bar{\mathbb{D}}}\psi(\bar{p}+r\,q)\,\mathrm{d}\lambda_{\bar{\mathbb{D}}}(q),
\end{equation}%
where $\bar{\mathbb{D}}$ is the closed unit disc centered at the origin and $\lambda_{\bar{\mathbb{D}}}$ denotes the standard volume measure on $\bar{\mathbb{D}}$. The value of $\Psi_r$ at $\bar{p}\in{P}$ is the average of $\psi$ over the disk of radius $r$ centered at the point $\bar{p}$. In the limit $r\to0$, the averaged function $\Psi_r$ converges locally uniformly to $\psi$ on $P$. In the following, the positive real number $r$ is the distance of the sensor to center of the unicycle. Note that the definition of $\Psi_r$ involves an integral over the closed unit disc $\bar{\mathbb{D}}$ while the definition of $\bar{G}_r$ involves an integral over the boundary of $\bar{\mathbb{D}}$, where $o\in\mathbb{S}$ is an outward-pointing unit vector. An application of the \emph{divergence theorem} (also known as \emph{Gauss's theorem} or \emph{Ostrogradsky's theorem}) leads to%
\begin{equation}\label{eq:14}
\bar{G}_r(\bar{p}) \ = \ \nabla\Psi_r(\bar{p})
\end{equation}%
for every $\bar{p}\in{P}$. Thus, on $P$, the averaged equation \cref{eq:12} can be written as the gradient system%
\begin{equation}\label{eq:15}
\dot{\bar{p}}(t) \ = \ -\varepsilon\,\tfrac{\lambda\,r}{2\,k}\,\nabla\Psi_r(\bar{p}(t)) + \varepsilon\,\bar{d}(t).
\end{equation}%
One can expect that the averaged signal function $\Psi_r$ has fewer undesired critical points than the original signal function $\psi$. In this sense, assumptions on $\Psi_r$ are weaker than assumptions on $\psi$. On the other hand, the shape of the function $\Psi_r$ does not only depend on $\psi$ but also on the distance $r$. For this reason, any assumption on $\Psi_r$ is an assumption on both $\psi$ and $r$. Moreover, in the context of source seeking, the signal function $\psi$ is an unknown quantity. Consequently, the user will not be able to check in advance whether $\psi$ has undesired local extrema. Therefore, the user can neither check conditions on $\psi$ nor conditions $\Psi_r$ prior to an implementation. In practice, assumptions on $\psi$ are as difficult to verify as assumptions on $\Psi_r$. We suppose that the averaged signal function $\Psi_r\colon{P}\to\mathbb{R}$ in \cref{eq:13} has the following property.%
\begin{assumption}\label{Assumption:1}
There exist real numbers $y_\ast$ and $y_0$ with $y_\ast<y_0$ such that the $y_0$-sublevel set of $\Psi_r$ is compact and such that the gradient of $\Psi_r$ is non-zero at every point $p\in{P}$ with $y_\ast<\Psi_r(p)\leq{y_0}$.
\end{assumption}%
In \Cref{Assumption:1}, the $y_0$-sublevel set of $\Psi_r$ will act as a set of admissible initial positions and the $y_\ast$-sublevel represents a set of target positions. If \Cref{Assumption:1} is satisfied, then the gradient system \cref{eq:15} has the following robust asymptotic stability property.%
\begin{lemma}\label{Lemma:1}
Suppose that \Cref{Assumption:1} is satisfied with $y_\ast$ and $y_0$ as therein. Then, there exist an open neighborhood $U$ of the $y_0$-sublevel set of $\Psi_r$ in $P$, a real number $y_0'>y_0$, a function $\beta$ of class $\mathcal{KL}$ with $\beta(s,0)=s$ for every $s\geq0$, a function $\gamma$ of class $\mathcal{K}$, and a positive real number $\delta$ such that the $y_0'$-sublevel set of the restriction of $\Psi_r$ to $U$ is compact and such that, for every $t_0\in\mathbb{R}$, every $\bar{p}_0\in{U}$ with $\Psi_r(\bar{p}_0)\leq{y_0'}$, every $\bar{d}\in{L_\infty^2}$ with $\|\bar{d}\|\leq\delta$, and every $\varepsilon>0$, the maximal solution $\bar{p}$ of \cref{eq:15} with initial value $\bar{p}_0$ at initial time $t_0$ satisfies%
\begin{subequations}\label{eq:16}%
\begin{align}
& \Psi_r(\bar{p}(t)) - y_\ast \ \leq \ \allowdisplaybreaks \\
& \quad \max\big\{\beta\big(\max\{\Psi_r(\bar{p}_0)-y_\ast,0\},\varepsilon\,(t-t_0)\big), \gamma(\|\bar{d}\|)\big\}
\end{align}%
\end{subequations}%
for every $t\geq{t_0}$.
\end{lemma}
The existence of $\beta\in\mathcal{KL}$ and $\gamma\in\mathcal{K}$ in \Cref{Lemma:1} follow from a similar reasoning as in the proof of Theorem~1 in \cite{Sontag1989}. In particular, one can interpret the property described by inequality \cref{eq:16} as a weak form of input-to-state stability with respect to sufficiently weak disturbances.

We end this section by establishing a direct link between the method in this paper and the approach in \cite{Tan2006,Tan2009} for steady-state extremum seeking in the presence of local extrema.%
\begin{remark}\label{Remark:1}
An averaged objective function also appears implicitly  in \cite{Tan2006,Tan2009}. The averaged ``reduced system'' in \cite{Tan2006,Tan2009} is a scalar system of the form%
\begin{equation}\label{eq:17}
\dot{\bar{\theta}}(t) \ = \ -\tfrac{\delta\,a}{2}\,\bar{g}_a(\bar{\theta}(t)),
\end{equation}%
where $a$ and $\delta$ are positive control parameters, the one-dimensional vector field $\bar{g}_a\colon\mathbb{R}\to\mathbb{R}$ is defined by%
\begin{equation}\label{eq:18}
\bar{g}_a(\bar{\theta}) \ := \ \tfrac{1}{\pi}\int_0^{2\pi}\tfrac{1}{a}\,\phi(\bar{\theta}+a\,\sin(\tau))\,\sin(\tau)\,\mathrm{d}\tau,
\end{equation}%
and $\phi\colon\mathbb{R}\to\mathbb{R}$ is a smooth steady-state input-output function. Substitution and integration by parts lead to%
\begin{equation}\label{eq:19}
\bar{g}_a(\bar{\theta}) \ = \ \nabla\Phi_a(\bar{\theta}),
\end{equation}%
where the averaged steady-state input-output function $\Phi_a$ is defined by%
\begin{equation}\label{eq:20}
\Phi_a(\bar{\theta}) \ := \ \tfrac{2}{\pi}\int_{-1}^{1}\sqrt{1-s^2}\,\phi(\bar{\theta}+a\,s)\,\mathrm{d}s.
\end{equation}%
This means that the averaged ``reduced system'' is driven into the negative gradient direction of an averaged steady-state input-output function. If the parameter $a$ is sufficiently large, then one can expect that $\Phi_a$ has fewer undesired critical points than $\phi$. The control parameter $a$ in \cite{Tan2006,Tan2009} corresponds to the sensor distance $r$ in our paper. We emphasize that this interpretation of the averaged ``reduced system'' does not appear in \cite{Tan2006,Tan2009}.
\end{remark}%


\section{Main result}\label{Section:4}
As indicated in \Cref{Section:3}, if the positive control parameter $\varepsilon$ is sufficiently small, then the position of the source-seeking unicycle \cref{eq:01}, \cref{eq:02} under control law \cref{eq:03}, \cref{eq:04} follows approximately the solutions of the gradient system \cref{eq:15} for the averaged signal function $\Psi_r\colon{P}\to\mathbb{R}$ in \cref{eq:13}. Moreover, under mild assumptions on $\Psi_r$, the gradient system has the robust asymptotic stability property in \Cref{Lemma:1}. This in turn will allow us to conclude a robust practical stability property for the approximating closed-loop system for sufficiently small $\varepsilon>0$. First, we introduce the following notation. For every real number $\eta$ and every nonempty set $\Lambda$ of real numbers, let $|\eta|_\Lambda$ denote the distance of $\eta$ to $\Lambda$; i.e., the infimum of all $|\eta-\tilde{\eta}|$ with $\tilde{\eta}\in\Lambda$.%
\begin{theorem}\label{Theorem:1}
Suppose that \Cref{Assumption:1} is satisfied with $y_\ast$ and $y_0$ as therein. Let $\Lambda$ be the set of $\psi(p+r\,o)$ with $p$ in the $y_0$-sublevel set of $\Psi_r$ and $o\in\mathbb{S}$. Let $\rho_v$, $\rho_\omega$, $\rho_\eta$ be positive real numbers. Then, there exists a function $\beta$ of class $\mathcal{KL}$ with $\beta(s,0)=s$ for every $s\geq0$, there exist functions $\gamma_{\text{r}}$, $\gamma_{\text{s}}$, $\gamma_{\text{t}}$, $\nu_p$, $\nu_\eta$ of class $\mathcal{K}$, and there exist positive constants $\varepsilon_0$, $\delta_{\text{r}}$, $\delta_{\text{s}}$, $\delta_{\text{t}}$, $\nu_v$ such that, for every $\varepsilon\in(0,\varepsilon_0)$, every $t_0\in\mathbb{R}$, every $p_0$ in the $y_0$-sublevel set of $\Psi_r$, every $o_0\in\mathbb{S}$, every $v_0\in\mathbb{R}$ with $|v_0|\leq\rho_v$, every $\omega_0\in\mathbb{R}$ with $|\omega_0-\omega_\ast|\leq\rho_\omega$, every $\eta_0\in\mathbb{R}$ with $|\eta_0|_{\Lambda}\leq\rho_\eta$, and all $d_{\text{r}},d_{\text{s}},d_{\text{t}}\in{L_\infty}$ with $\|d_{\text{r}}\|\leq\delta_{\text{r}}$, $\|d_{\text{s}}\|\leq\delta_{\text{s}}$, $\|d_{\text{t}}\|\leq\delta_{\text{t}}$, the maximal solution $(p,o,v,\omega,\eta)$ of the closed-loop system \cref{eq:06} with initial value $(p_0,o_0,v_0,\omega_0,\eta_0)$ at initial time $t_0$ satisfies%
\begin{subequations}\label{eq:21}%
\begin{align}
\Psi_r(p(t)) - y_\ast & \ \leq \ \label{eq:21:a}  \allowdisplaybreaks \\
& \!\!\!\!\!\!\!\!\!\!\!\!\!\!\!\!\!\!\!\!\!\!\!\!\max\big\{\beta\big(\max\{\Psi_r(p_0)-y_\ast,0\},\varepsilon\,(t-t_0)\big), \label{eq:21:b} \allowdisplaybreaks \\
& \!\!\!\!\!\!\!\!\!\!\!\!\!\!\!\!\!\!\!\!\!\!\!\!\!\hphantom{\max{m}} \gamma_{\text{r}}(\|d_{\text{r}}\|) + \gamma_{\text{s}}(\|d_{\text{s}}\|) + \gamma_{\text{t}}(\|d_{\text{t}}\|)\big\} + \nu_p(\varepsilon), \label{eq:21:c} \allowdisplaybreaks \\
|v(t)| & \ \leq \ |v_0|\,\mathrm{e}^{-\frac{1}{\varepsilon}\,\frac{k}{m}\,(t-t_0)} + \nu_v\,\varepsilon, \label{eq:21:d} \allowdisplaybreaks \\
|\omega(t)-\omega_\ast| & \ \leq \ |\omega_0-\omega_\ast|\,\mathrm{e}^{-\frac{\kappa}{J}\,(t-t_0)} + \tfrac{1}{\kappa}\,\|d_{\text{r}}\|, \label{eq:21:e} \allowdisplaybreaks \\
|\eta(t)|_\Lambda & \ \leq \ |\eta_0|_\Lambda\,\mathrm{e}^{-\varepsilon\,h\,(t-t_0)} + \|d_{\text{s}}\| + \nu_\eta(\varepsilon) \label{eq:21:f}
\end{align}%
\end{subequations}%
for every $t\geq{t_0}$.
\end{theorem}%
The proof of \Cref{Theorem:1} in \hyperlink{appendixChapter}{Appendix} will reveal that the function $\beta$ in \cref{eq:21:c} is the same as in \Cref{Lemma:1} and that the functions $\gamma_{\text{r}},\gamma_{\text{s}},\gamma_{\text{t}}$ in \cref{eq:21:c} originate from the function $\gamma$ in \Cref{Lemma:1}. Therefore, one may say that certain robustness and stability properties of the averaged gradient system \cref{eq:15} carry over to the approximating closed-loop system \cref{eq:06}. Good robustness and stability properties can be expected if the the averaged signal function has a strong gradient vector field. An increase of the sensor distance $r$ can be helpful to strengthen the gradient of the averaged signal function in the presence of local extrema.


%
\begin{figure}[t]%
\centering$\quad$\includegraphics{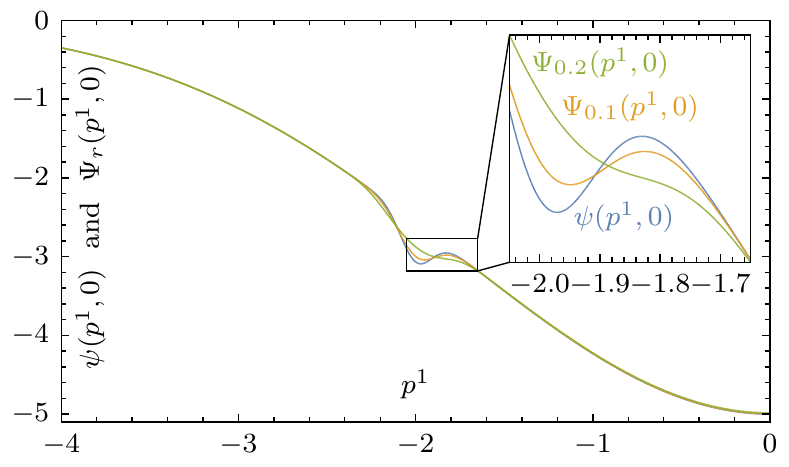}%
\caption{Graphs of the sensed signal function $\psi$ (blue) given by \cref{eq:22} and the averaged signal function $\Psi_r$ defined in \cref{eq:13} for sensor distances $r=0.1$ (orange) and $r=0.2$ (green) to the vehicle center. To allow a visual comparison with \Cref{Figure:3}, the functions are shown for arguments $p=(p^1,p^2)\in\mathbb{R}^2$ with $p^1\in[-4,0]$ and $p^2=0$.}%
\label{Figure:2}%
\end{figure}%
\begin{figure}[t]%
\centering$\begin{matrix}\!\!\!\!\includegraphics{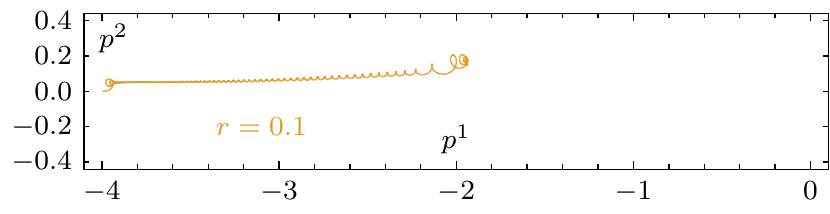}\\\!\!\!\!\includegraphics{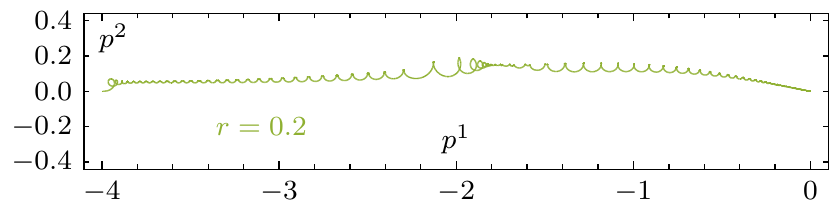}\end{matrix}$%
\caption{Trajectory of the unicycle's center position $p$ for two different distances $r$ of the sensor to $p$. In each of the simulations, the initial position is the point $(-4,0)$ and the source is at the origin. The choice of colors is the same as in \Cref{Figure:2}.}%
\label{Figure:3}%
\end{figure}%
\section{Numerical test}\label{Section:5}
We test the proposed method for the case that the measured signal function $\psi$ is given by%
\begin{equation}\label{eq:22}
\psi(p) \ = \ - 5\,\mathrm{e}^{-|p|^2/6} - \tfrac{1}{2}\,\mathrm{e}^{-4 (|p|^2-2^2)^2},
\end{equation}%
where $|\cdot|$ is the Euclidean norm on $\mathbb{R}^2$. The graph of this radially symmetric function is indicated in \Cref{Figure:2}. One can see that $\psi$ attains its global minimum value at the origin, but there is also a circular set of local minimizers with distance $\approx{2}$ to the origin. If we would use a method that steers the unicycle into the negative gradient direction of $\psi$, then the circular set of local minimizers would be an insurmountable barrier. In this case an initial position with distance $>2$ to the origin, would only lead to convergence into the circular set of local minimizers, but not to convergence toward the source at the origin. The proposed control law \cref{eq:03}, \cref{eq:04} has the ability to pass through the blocking set of local minimizers because it steers the unicycle into the negative gradient direction of the averaged signal function $\Psi_r$ in \cref{eq:13}. Note that $\Psi_r$ depends on the distance $r>0$ of the sensor to the center of the vehicle. For this reason, the shape of $\Psi_r$ depends on the distance $r$. The larger $r$, the larger is the area of averaging in \cref{eq:13}. It is natural to expect that the circular set of local minima of $\psi$ is eliminated if the area of averaging is sufficiently large. This is exactly what we can observe in \Cref{Figure:2}. For $r=0.1$, the averaging effect leads to a less distinct circular set of local minima but the averaging effect is too weak to eliminate the local minimizers completely. For $r=0.2$, the area of averaging is large enough and the corresponding averaged signal function only has one minimizer at the origin. Thus, we expect that the unicycle cannot pass through the circular set of local minima of $\psi$ for small values of the control parameter $\varepsilon$ if $r=0.1$. On the other hand, \Cref{Theorem:1} guarantees convergence of the position into a small neighborhood of the origin for sufficiently small values of the control parameter $\varepsilon$ if $r=0.2$.

\begin{figure}[t]%
\centering\includegraphics{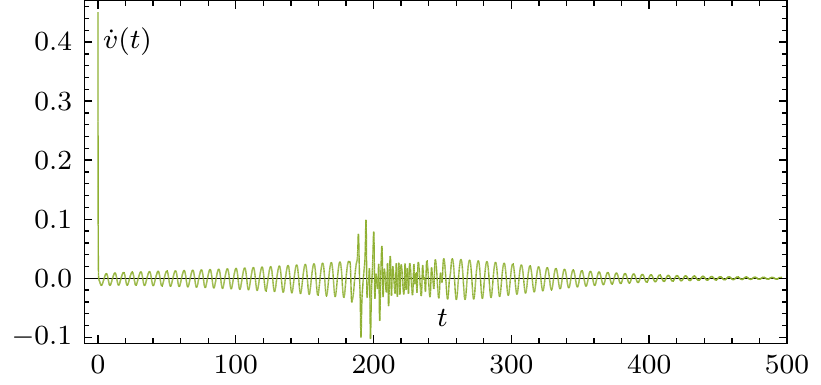}%
\caption{Translational acceleration $t\mapsto\dot{v}(t)$ for $r=0.2$.}%
\label{Figure:4}%
\end{figure}

\begin{figure}[t]%
\centering\includegraphics{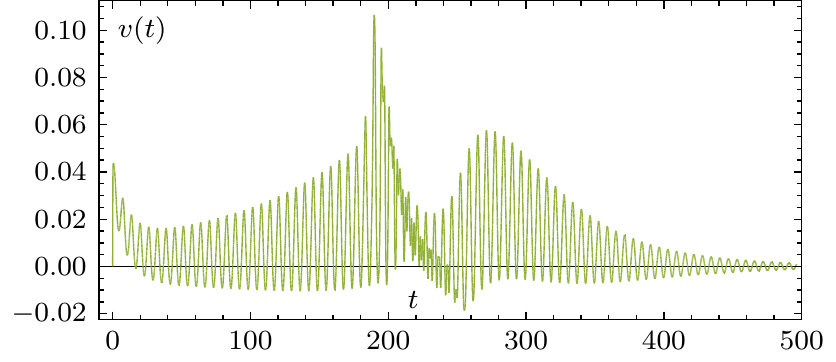}%
\caption{Translational velocity $t\mapsto{v(t)}$ for $r=0.2$.}%
\label{fig:v}%
\end{figure}

To generate numerical data, we set all constants, except for $r$, equal to $1$, and we choose the control parameter $\varepsilon$ equal to $0.1$. The simulation results in \Cref{Figure:3} indicate the position of the unicycle under control law \cref{eq:03}, \cref{eq:04} with initial condition $p(0)=(-4,0)$, $o(0)=(1,0)$, and $v(0)=\omega(0)=\eta(0)=0$. One can observe in \Cref{Figure:3} that the unicycle gets stuck in the circular set of local minima if the distance $r$ of the sensor to vehicle center is equal to $0.1$. If $r$ is increased to $0.2$, then, as expected in the previous paragraph, the unicycle reaches the source at the origin. \Cref{Figure:4,fig:v} also show the translational acceleration and the translational velocity for sensor distance distance $r=0.2$. One can see in \Cref{fig:v} that the translational velocity has a positive average. This can be explained as follows. Since our control law steers the unicycle towards a minimum of the signal function, the sensed signal value $\hat{y}(t)$ decreases and the high-pass filter state $\eta(t)$ follows slowly. This implies that the difference $-\lambda(\hat{y}(t)-\eta(t))$ in \cref{eq:03:a} is, on average, positive, which in turn causes a translational velocity with positive average. One can also observe in \Cref{fig:v} that the average value is close to zero for time parameters between 200 to 250. In this time interval, the unicycle is close to the local minima of the signal function and therefore the gradient of the averaged signal function is very weak.%


\appendix
\section*{Appendix: Proof of \texorpdfstring{\Cref{Theorem:1}}{Theorem~\ref{Theorem:1}}}\hypertarget{appendixChapter}{}

The proof of \Cref{Theorem:1} consists of the following two steps. In the first step, we perform a first-order averaging transformation as, e.g., in Section~9.1 of the textbook \cite{BulloBook} to prove that the position variable in \cref{eq:06:a} of the closed-loop system approximates the solutions of the averaged system \cref{eq:12} with decreasing parameter $\varepsilon>0$. This approximation property is made precise in \Cref{Proposition:1} below. In the second step, we use \Cref{Proposition:1} to prove the asserted robustness and stability properties of the closed-loop system. Recall from \Cref{Section:3} that the averaged system \cref{eq:12} is the negative gradient system \cref{eq:15} for the averaged signal function. Moreover, if \Cref{Assumption:1} is satisfied, then the negative gradient system has the robustness and stability properties in \Cref{Lemma:1}. Using the approximation property from the first step, we then show that the robustness and stability properties of the averaged system carry over to the approximating closed-loop system if the parameter $\varepsilon>0$ is sufficiently small. The rather technical proof of robustness and stability for the closed-loop system is split up into a sequence of auxiliary results (\Cref{Lemma:2,Lemma:3,Lemma:4}).

As explained above, we first prove the following approximation result for the position of the unicycle.%
\begin{proposition}\label{Proposition:1}
Let $K$ be a compact subset of $M$ and let $\delta_{\text{r}},\delta_{\text{s}},\delta_{\text{t}}>0$. Then there exist $c_{\text{r}},c_{\text{s}},c_{\text{t}},c_1,c_2,\varepsilon_0>0$ such that, for every $\varepsilon\in(0,\varepsilon_0)$, all $d_{\text{r}},d_{\text{s}},d_{\text{t}}\in{L_\infty}$ with $\|d_{\text{r}}\|\leq\delta_{\text{r}}$, $\|d_{\text{s}}\|\leq\delta_{\text{s}}$, $\|d_{\text{t}}\|\leq\delta_{\text{t}}$, every solution $x=(p,o,v,\omega,\eta)\colon{I}\to{M}$ of \cref{eq:06}, and all $t_1,t_3\in{I}$ with $t_1<t_3$, the following implication holds: If $x(t)\in{K}$ for every $t\in[t_1,t_3]$, then there exists an absolutely continuous map $\bar{d}\colon[t_1,t_3]\to\mathbb{R}^2$ such that%
\begin{equation}\label{eq:23}
|\bar{d}(t)| \ \leq \ c_{\text{r}}\,\|d_{\text{r}}\| + c_{\text{s}}\,\|d_{\text{s}}\| + c_{\text{t}}\,\|d_{\text{t}}\|
\end{equation}%
for every $t\in[t_1,t_3]$ and%
\begin{subequations}\label{eq:24}%
\begin{align}
& \Big|p(t_2)-p(t_1) - \int_{t_1}^{t_2}\big( - \varepsilon\,\tfrac{\lambda\,r}{2\,k}\bar{G}_r(p(t)) + \varepsilon\,\bar{d}(t)\big)\,\mathrm{d}t\Big| \allowdisplaybreaks \\
& \quad \leq \ \varepsilon\,c_1 + \varepsilon^2\,c_2\,|t_2-t_1|
\end{align}%
\end{subequations}%
for every $t_2\in[t_1,t_3]$.
\end{proposition}
\begin{proof}
Choose sufficiently large $\rho_p,\rho_v,\rho_\omega,\rho_\eta>1$ such that $|p|\leq\rho_p-1$ $|v|\leq\rho_v-1$, $|\omega-\omega_\ast|\leq\rho_\omega$, $|\eta|\leq\rho_\eta$ for every $(p,o,v,\omega,\eta)\in{K}$. Let $\tilde{K}$ be the compact set of $(\tilde{p},\tilde{o},\tilde{v},\tilde{\Omega},\tilde{\eta})$ in $M$ with $\tilde{p}\in\rho_p\bar{\mathbb{D}}$, $\tilde{o}\in\mathbb{S}$, $|\tilde{v}|\leq\rho_v$, $|\tilde{\omega}|\leq\rho_\omega$, $|\tilde{\eta}|\leq\rho_\eta$. Recall that the maps $\tilde{G}_r$ and $\bar{G}_r$ are defined in \cref{eq:09} and \cref{eq:11}, respectively. For every $t\in\mathbb{R}$, every $\tilde{o}\in\mathbb{S}$, and every $\tilde{\eta}\in\mathbb{R}$, define a vector field $H[t,\tilde{o},\tilde{\eta}]$ on $\mathbb{R}^2$ by%
\begin{equation*}
H[t,\tilde{o},\tilde{\eta}](\tilde{p}) \ := \ -\tfrac{\lambda\,r}{2\,k}\int_0^t\big(\tilde{G}_r(s,\tilde{p},\tilde{o},\tilde{\eta}) - \bar{G}_r(\tilde{p})\big)\,\mathrm{d}s    
\end{equation*}%
Because of \cref{eq:10} and the $(2\pi/\omega_\ast)$-periodicity of $\tilde{G}_r$ in the first argument, also $t\mapsto{H[t,\tilde{o},\tilde{\eta}](\tilde{p})}$ is $(2\pi/\omega_\ast)$-periodic. Let $(\tau,\tilde{p})\mapsto\Phi^{H[t,\tilde{o},\tilde{\eta}]}_\tau(\tilde{p})$ denote the flow of $H[t,\tilde{o},\tilde{\eta}]$. Then, there is a sufficiently small $\varepsilon_0>0$ such that $|\Phi^{H[t,\tilde{o},\tilde{\eta}]}_{-\varepsilon}(p)|\leq\rho_p$ and $\varepsilon\frac{\lambda}{k}|\psi(p+r\tilde{o})-\tilde{\eta}|\leq{1}$ for every $\varepsilon\in[0,\varepsilon_0]$, every $p\in(\rho_p-1)\bar{\mathbb{D}}$, every $\tilde{o}\in\mathbb{S}$, and every $\tilde{\eta}\in[-\rho_\eta,\rho_\eta]$. Furthermore, for every $\varepsilon\in(0,\varepsilon_0)$, all $d_{\text{r}},d_{\text{s}},d_{\text{t}}\in{L_\infty}$, and every solution $x=(p,o,v,\omega,\eta)\colon{I}\to{M}$ of \cref{eq:06}, the following implication holds: If $x(t)\in{K}$ for every $t\in{I}$, then a well-defined curve $\tilde{x}=(\tilde{p},\tilde{o},\tilde{v},\tilde{\omega},\tilde{\eta})\colon{I}\to{M}$ with values in $\tilde{K}$ is defined by%
\begin{subequations}\label{eq:25}%
\begin{align}
\tilde{p}(t) & \ := \ \Phi^{H[t,\tilde{o}(t),\tilde{\eta}(t)]}_{-\varepsilon}(p(t)), \label{eq:25:a} \allowdisplaybreaks \\
\tilde{o}(t) & \ := \ \exp(-\omega_\ast\,t\,R)\,o(t), \allowdisplaybreaks \\
\tilde{v}(t) & \ := \ v(t) + \varepsilon\,\tfrac{\lambda}{k}\,\big(\psi(p(t)+r\,o(t))-\eta(t)\big), \allowdisplaybreaks \\
\tilde{\omega}(t) & \ := \ \omega(t) - \omega_\ast, \allowdisplaybreaks \\
\tilde{\eta}(t) & \ := \ \eta(t).
\end{align}%
\end{subequations}%
Note that \cref{eq:25} is almost the same change of variables as \cref{eq:07} except for the additional transformation in \cref{eq:25:a}, which will allow us to prove the asserted approximation of \cref{eq:12}. Some lengthy computations using known rules for the derivatives of flows (see, e.g., Exercise~E9.4 and Proposition~3.85 in \cite{BulloBook}) lead to the result that there exist smooth maps $A,B\colon[0,\varepsilon_0]\times\mathbb{R}\times\tilde{K}\to\mathbb{R}^2$ and $C,D\colon[0,\varepsilon_0]\times\mathbb{R}\times\tilde{K}\times\mathbb{R}^3\to\mathbb{R}$, which are $(2\pi)/\omega_\ast$-periodic in the second argument, such that, for every $\varepsilon\in(0,\varepsilon_0)$, every $d=(d_{\text{r}},d_{\text{s}},d_{\text{t}})\in{L_\infty^3}$, and every solution $x=(p,o,v,\omega,\eta)\colon{I}\to{M}$ of \cref{eq:06}, the following implication holds: If $x(t)\in{K}$ for every $t\in{I}$, then the components $\tilde{p},\tilde{v},\tilde{\omega}$ of $\tilde{x}=(\tilde{p},\tilde{o},\tilde{v},\tilde{\omega},\tilde{\eta})\colon{I}\to{M}$ defined by \cref{eq:25} satisfy the differential equations%
\begin{align*}
\dot{\tilde{p}}(t) & \ = \ -\varepsilon\,\tfrac{\lambda\,r}{2\,k}\,\bar{G}_r(\tilde{p}(t)) + \tilde{v}(t)\,A(\varepsilon,t,\tilde{x}(t)) \allowdisplaybreaks \\
& \qquad + \varepsilon\,\tilde{\omega}(t)\,B(\varepsilon,t,\tilde{x}(t)) + \varepsilon^2\,C(\varepsilon,t,\tilde{x}(t),d(t)), \allowdisplaybreaks \\
\dot{\tilde{v}}(t) & \ = \ -\tfrac{1}{\varepsilon}\,\tfrac{k}{m}\,\tilde{v}(t) - \tfrac{\lambda}{m}\,d_{\text{s}}(t) + \tfrac{1}{m}\,d_{\text{t}}(t) + \varepsilon\,D(\varepsilon,t,\tilde{x}(t),d(t)), \allowdisplaybreaks \\
\dot{\tilde{\omega}}(t) & \ = \ -\tfrac{\kappa}{J}\,\tilde{\omega}(t) + \tfrac{1}{J}\,d_{\text{r}}(t).
\end{align*}%
Since $A$, $B$, $C$, $D$, and $(\varepsilon,t,p,\tilde{o},\tilde{\eta})\mapsto\Phi^{H[t,\tilde{o},\tilde{\eta}]}_{-\varepsilon}(p)$ are smooth and periodic in the second argument, there exist sufficiently large $c_A,c_B,c_C,c_D,c_H>0$ such that%
\begin{align*}
& |A(\varepsilon,t,\tilde{x})| \, \leq \, c_A, \ \ |B(\varepsilon,t,\tilde{x})| \, \leq \, c_B, \ \ |C(\varepsilon,t,\tilde{x},d)| \, \leq \, c_C,  \allowdisplaybreaks \\
& |D(\varepsilon,t,\tilde{x},d)| \, \leq \, c_D, \ \ |p-\Phi^{H[t,\tilde{o},\tilde{\eta}]}_{-\varepsilon}(p)| \leq \varepsilon\,c_H
\end{align*}%
for every $\varepsilon\in[0,\varepsilon_0]$, every $t\in\mathbb{R}$, every $p\in(\rho_p-1)\bar{\mathbb{D}}$, every $\tilde{x}=(\tilde{p},\tilde{o},\tilde{v},\tilde{\omega},\tilde{\eta})\in\tilde{K}$, and every $d=(d_{\text{r}},d_{\text{s}},d_{\text{t}})\in\mathbb{R}^3$ with $|d_{\text{r}}|\leq\delta_{\text{r}}$, $|d_{\text{s}}|\leq\delta_{\text{s}}$, $|d_{\text{t}}|\leq\delta_{\text{t}}$. Moreover, since $\bar{G}_r$ is smooth, there is a Lipschitz constant $L_G$ for $\bar{G}_r$ on $\rho_p\bar{\mathbb{D}}$. Finally,  set%
\begin{align*}
& c_{\text{r}} \, := \, \tfrac{1}{\kappa}\,c_B, \quad c_{\text{s}} := \tfrac{\lambda}{k}\,c_A, \quad c_{\text{t}} \, := \, \tfrac{1}{k}\,c_A, \allowdisplaybreaks \\
& c_1 \, := \, \tfrac{m}{k}\,\rho_v\,c_A + \tfrac{J}{\kappa}\,\rho_\omega\,c_B + 2\,c_H, \allowdisplaybreaks \\
& c_2 \, := \, \tfrac{m}{k}\,c_A\,c_D + c_C + \tfrac{\lambda\,r}{2\,k}\,L_G\,c_H.
\end{align*}%
Next, we indicate that the asserted implication in \Cref{Proposition:1} holds for above choice of $c_{\text{r}}$, $c_{\text{s}}$, $c_{\text{t}}$, $c_1$, $c_2$, $\varepsilon_0$.

Fix an arbitrary $\varepsilon\in(0,\varepsilon_0)$, arbitrary $d_{\text{r}},d_{\text{s}},d_{\text{t}}\in{L_\infty}$ with $\|d_{\text{r}}\|\leq\delta_{\text{r}}$, $\|d_{\text{s}}\|\leq\delta_{\text{s}}$, $\|d_{\text{t}}\|\leq\delta_{\text{t}}$, an arbitrary solution $x=(p,o,v,\omega,\eta)\colon{I}\to{M}$ of \cref{eq:06}, and arbitrary $t_1,t_3\in{I}$ with $t_1<t_3$. Assume that $x(t)\in{K}$ for every $t\in[t_1,t_3]$. Define $\tilde{x}=(\tilde{p},\tilde{o},\tilde{v},\tilde{\omega},\tilde{\eta})\colon[t_1,t_3]\to{M}$ by \cref{eq:25}. Then an absolutely continuous map $\bar{d}\colon[t_1,t_3]\to\mathbb{R}^2$ is defined by%
\begin{align*}
& \bar{d}(t) \ := \ \int_{t_1}^td_{\text{r}}(s)\,\tfrac{1}{J}\,\mathrm{e}^{-\frac{\kappa}{J}\,(t-s)}\,\mathrm{d}s\,B(\varepsilon,t,\tilde{x}(t)) \allowdisplaybreaks \\
& \ + \int_{t_1}^t\big( - \lambda\,d_{\text{s}}(s)+d_{\text{t}}(s)\big)\tfrac{1}{\varepsilon}\,\tfrac{1}{m}\,\mathrm{e}^{-\frac{1}{\varepsilon}\,\frac{k}{m}\,(t-s)}\,\mathrm{d}s\,A(\varepsilon,t,\tilde{x}(t)).
\end{align*}%
It is obvious that $\bar{d}$ satisfies \cref{eq:23} for every $t\in[t_1,t_3]$. Moreover, using the variation of constants formula to compute $\tilde{v}$ and $\tilde{\omega}$, one can verify that \cref{eq:24} is satisfied for every $t_2\in[t_1,t_3]$
\end{proof}
Now we are ready to prove \Cref{Theorem:1}. We will use the approximation result, \Cref{Proposition:1}, to show that asymptotic stability for the averaged gradient system \cref{eq:15} implies practical asymptotic stability for the approximating closed-loop system \cref{eq:06}. From now on, we suppose that \Cref{Assumption:1} is satisfied with $y_\ast$ and $y_0$ as therein. Then, there exist $U$, $y_0'$, $\beta$, $\gamma$, and $\delta$ as in \Cref{Lemma:1}. Following the notation in \Cref{Theorem:1}, we let $\Lambda$ denote the set of $\psi(p+r\,o)$ with $p$ in the $y_0$-sublevel set of $\Psi_r$ and $o\in\mathbb{S}$. Let $\rho_v$, $\rho_\omega$, $\rho_\eta$ be positive real numbers. We have to show the existence of class-$\mathcal{K}$ functions $\gamma_{\text{r}}$, $\gamma_{\text{s}}$, $\gamma_{\text{t}}$, $\nu_p$, $\nu_\eta$ and positive constants $\varepsilon_0$, $\delta_{\text{r}}$, $\delta_{\text{s}}$, $\delta_{\text{t}}$, $\nu_v$ as in \Cref{Theorem:1}.

Let $C_0$ be the compact set of $p\in{P}$ with $\Psi_r(p)\leq{y_0}$. Let $C$ be the compact set of $p\in{U}$ with $\Psi_r(p)\leq{y_0'}$. Let $\Lambda'$ be the compact set of $\psi(p+ro)$ with $p\in{C}$ and $o\in\mathbb{S}$. Let $K_0$ be the compact set of $(p_0,o_0,v_0,\Omega_0,\eta_0)$ in $M$ with $p_0\in{C_0}$, $o_0\in\mathbb{S}$, $|v_0|\leq\rho_v$, $|\omega_0-\omega_\ast|\leq\rho_\omega$, $|\eta_0|_\Lambda\leq\rho_\eta$. Let $K$ be the compact set of $(p,o,v,\Omega,\eta)$ in $M$ wit $p\in{C}$, $o\in\mathbb{S}$, $|v|\leq\rho_v+1$, $|\omega-\omega_\ast|\leq\rho_\omega+1$, $|\eta|_{\Lambda'}\leq\rho_\eta+1$. Using the variation of constants formula for the linear subsystems \cref{eq:06:c,eq:06:d,eq:06:e,eq:06:f} one can easily verify the following statement.%
\begin{lemma}\label{Lemma:2}
There exist  $\varepsilon_0,\delta_{\text{r}}, \delta_{\text{s}},\delta_{\text{t}},\nu_v >0$ such that, for every $\varepsilon\in(0,\varepsilon_0)$, all $t_0,t_1\in\mathbb{R}$ with $t_0\leq{t_1}$, every $x_0=(p_0,o_0,v_0,\omega_0,\eta_0)$ in $K_0$ and all $d_{\text{r}},d_{\text{s}},d_{\text{t}}\in{L_\infty}$ with $\|d_{\text{r}}\|\leq\delta_{\text{r}}$, $\|d_{\text{s}}\|\leq\delta_{\text{s}}$, $\|d_{\text{t}}\|\leq\delta_{\text{t}}$, the maximal solution $x=(p,o,v,\omega,\eta)\colon{I}\to{M}$ of \cref{eq:06} with $x(t_0)=x_0$ satisfies the following implication: If $[t_0,t_1]\subset{I}$ and if $p(t)\in{C}$ for every $t\in[t_0,t_1]$, then $x(t)\in{K}$ for every $t\in[t_0,t_1]$ and the estimates \cref{eq:21:d,eq:21:e} hold for every $t\in[t_0,t_1]$.
\end{lemma}%
For the rest of the proof of \Cref{Theorem:1}, select $\varepsilon_0$, $\delta_{\text{r}}$, $\delta_{\text{s}}$, $\delta_{\text{t}}$, $\nu_v$ as in \Cref{Lemma:2}. After possibly reducing $\varepsilon_0>0$, we can also select $c_{\text{r}},c_{\text{s}},c_{\text{t}},c_1,c_2>0$ as in \Cref{Proposition:1}. Recall that $\gamma$ and $\delta$ are given by \Cref{Lemma:1}. Equation~(12) in \cite{Sontag1989} implies the existence of class-$\mathcal{K}$ functions $\gamma_{\text{r}},\gamma_{\text{s}},\gamma_{\text{t}}$ such that%
\begin{equation*}
\gamma(c_{\text{r}}|d_{\text{r}}| + c_{\text{s}}|d_{\text{s}}| + c_{\text{t}}|d_{\text{t}}|) \ \leq \ \gamma_{\text{r}}(|d_{\text{r}}|) + \gamma_{\text{s}}(|d_{\text{s}}|) + \gamma_{\text{t}}(|d_{\text{t}}|)
\end{equation*}%
for all $d_{\text{r}},d_{\text{s}},d_{\text{t}}\in\mathbb{R}$. After possibly reducing $\delta_{\text{r}},\delta_{\text{s}},\delta_{\text{t}}>0$, we can ensure that $c_{\text{r}}\delta_{\text{r}}+c_{\text{s}}\delta_{\text{s}}+c_{\text{t}}\delta_{\text{t}} \leq \delta$ and that
\begin{equation*}
\gamma_{\text{r}}(\delta_{\text{r}}) + \gamma_{\text{s}}(\delta_{\text{s}}) + \gamma_{\text{t}}(\delta_{\text{t}}) \ \leq \ \tfrac{1}{3}\,\min\{y_0-y_\ast,y_0'-y_0\}.
\end{equation*}
Since $\Psi_r$ is smooth, there exist Lipschitz constants $L_\Psi$ and $L_G$ for $\Psi_r$ and $\bar{G}_r$ on the compact subset $C$ of $P$. After possibly reducing $\varepsilon_0>0$, we can ensure that $T_\varepsilon:=\frac{1}{\varepsilon}\log(\log\frac{1}{\varepsilon})$ is a well-defined positive real number for every $\varepsilon\in(0,\varepsilon_0)$. Using the definition of $T_\varepsilon$, one can verify that there exists a function $\nu_p$ of class $\mathcal{K}$ such that%
\begin{equation*}
\max\big\{L_\Psi\,\big(\varepsilon{c_1} + \varepsilon^2\,c_2\,T_\varepsilon\big)\,\mathrm{e}^{\varepsilon\,\frac{\lambda{r}}{2k}\,L_G\,T_\varepsilon},\beta(y_0-y_\ast,\varepsilon\,T_\varepsilon)\big\} \leq \tfrac{1}{3}\,\nu_p(\varepsilon)
\end{equation*}%
for every $\varepsilon\in(0,\varepsilon_0)$. After possibly reducing $\varepsilon_0>0$, we  ensure  $\nu_p(\varepsilon)\leq\min\{y_0'-y_0,y_0-y_\ast\}/2$ for every $\varepsilon\in(0,\varepsilon_0)$.%
\begin{lemma}\label{Lemma:3}
Let $\varepsilon\in(0,\varepsilon_0)$, let $t_0,t_1\in\mathbb{R}$ with $t_0\leq{t_1}$, let $x_0\in{K_0}$, let $d_{\text{r}},d_{\text{s}},d_{\text{t}}\in{L_\infty}$ with $\|d_{\text{r}}\|\leq\delta_{\text{r}}$, $\|d_{\text{s}}\|\leq\delta_{\text{s}}$, $\|d_{\text{t}}\|\leq\delta_{\text{t}}$, and let $x=(p,o,v,\omega,\eta)\colon{I}\to{M}$ be the maximal solution of the closed-loop system \cref{eq:06} with $x(t_0)=x_0$. Assume that $[t_0,t_1]\subset{I}$, that $x(t)\in{K}$ for every $t\in[t_0,t_1]$, and that $p(t_1)\in{C_0}$. Then $[t_1,t_1+T_\varepsilon]\subset{I}$ and $x(t_2)\in{K}$ and%
\begin{align*}
&\Psi(p(t_2))-y_\ast \ \leq \ \max\big\{\beta\big(\max\{\Psi_r(p(t_1))-y_\ast,0\}, \allowdisplaybreaks \\
& \qquad \varepsilon\,(t_2-t_1)\big),\,\gamma_{\text{r}}(\|d_{\text{r}}\|) + \gamma_{\text{s}}(\|d_{\text{s}}\|) + \gamma_{\text{t}}(\|d_{\text{t}}\|)\big\} + \tfrac{1}{3}\nu_p(\varepsilon)
\end{align*}%
for every $t_2\in[t_1,t_1+T_\varepsilon]$.
\end{lemma}%
\begin{proof}
Let $t_3$ be the largest element of $[t_1,t_1+T_\varepsilon]$ with $[t_1,t_3]\subset{I}$ and $p(t)\in{C}$ for every $t\in[t_1,t_3]$. Then, we obtain from \Cref{Lemma:2} that $x(t)\in{K}$ for every $t\in[t_1,t_3]$.  \Cref{Proposition:1} states that there is an absolutely continuous map $\bar{d}\colon[t_1,t_3]\to\mathbb{R}^2$ as therein. Let $\bar{p}$ be the maximal solution of \cref{eq:15} with $\bar{p}(t_1)=p(t_1)$. By \Cref{Lemma:1} and \cref{eq:23}, we have%
\begin{align*}
&\Psi(\bar{p}(t_2))-y_\ast \ \leq \ \max\big\{\beta\big(\max\{\Psi_r(p(t_1))-y_\ast,0\}, \allowdisplaybreaks \\
& \qquad \varepsilon\,(t_2-t_1)\big),\,\gamma_{\text{r}}(\|d_{\text{r}}\|) + \gamma_{\text{s}}(\|d_{\text{s}}\|) + \gamma_{\text{t}}(\|d_{\text{t}}\|)\big\}
\end{align*}%
for every $t_2\in[t_1,t_3]$. On the other hand, estimate \cref{eq:24} and the Gronwall inequality in integral form lead to%
\begin{equation*}
L_\Psi\,|p(t_2)-\bar{p}(t_2)| \, \leq \, L_\Psi\,\big(\varepsilon{c_1} + \varepsilon^2\,c_2\,T_\varepsilon\big)\,\mathrm{e}^{\varepsilon\,\frac{\lambda{r}}{2k}\,L_G\,T_\varepsilon} \, \leq \, \tfrac{1}{3}\,\nu_p(\varepsilon)
\end{equation*}%
for every $t_2\in[t_1,t_3]$. Together, this implies that the asserted estimate is true for every $t_2\in[t_1,t_3]$. In particular, we obtain that $p(t_2)$ is an element of the interior of $C$ for every $t_2\in[t_1,t_3]$ and therefore $t_3=t_1+T_\varepsilon$.
\end{proof}
Next, one can apply \Cref{Lemma:3} repeatedly on time intervals of length $T_\varepsilon$. A simple proof by induction leads to the following result.%
\begin{lemma}\label{Lemma:4}
Let $\varepsilon\in(0,\varepsilon_0)$, let $t_0\in\mathbb{R}$, let $x_0=(p_0,o_0,v_0,$ $\omega_0,\eta_0)$ in $K_0$, let $d_{\text{r}},d_{\text{s}},d_{\text{t}}\in{L_\infty}$ with $\|d_{\text{r}}\|\leq\delta_{\text{r}}$, $\|d_{\text{s}}\|\leq\delta_{\text{s}}$, $\|d_{\text{t}}\|\leq\delta_{\text{t}}$, and let $x=(p,o,v,\omega,\eta)\colon{I}\to{M}$ be the maximal solution of the closed-loop system \cref{eq:06} with $x(t_0)=x_0$. Then $[\vphantom{]}t_0,\infty\vphantom{(})\subset{I}$, $x(t)\in{K}$ for every $t\geq{t_0}$, and estimate \cref{eq:21:a,eq:21:b,eq:21:c} holds for every $t\geq{t_0}$.
\end{lemma}%
It follows from the continuity of $\psi$ and $\Psi_r$, and from the definition of $\Lambda$ that there exists a function $\nu_\eta$ of class $\mathcal{K}$ such that%
\begin{equation*}
\psi(p+ro) \ \leq \ \max\Lambda + \nu_\eta(\varepsilon)
\end{equation*}%
for every $\varepsilon\in(0,\varepsilon_0)$, every $p\in{U}$ with $\Psi_r(p)\leq{y_0+\nu_p(\varepsilon)}$, and every $o\in\mathbb{S}$. To complete the proof of \Cref{Theorem:1}, let $\varepsilon\in(0,\varepsilon_0)$, let $t_0\in\mathbb{R}$, let $x_0=(p_0,o_0,v_0,\omega_0,\eta_0)$ in $K_0$, let $d_{\text{r}},d_{\text{s}},d_{\text{t}}\in{L_\infty}$ with $\|d_{\text{r}}\|\leq\delta_{\text{r}}$, $\|d_{\text{s}}\|\leq\delta_{\text{s}}$, $\|d_{\text{t}}\|\leq\delta_{\text{t}}$, and let $x=(p,o,v,\omega,\eta)\colon{I}\to{M}$ be the maximal solution of the closed-loop system \cref{eq:06} with $x(t_0)=x_0$. Then, by \Cref{Lemma:4}, $[\vphantom{]}t_0,\infty\vphantom{(})\subset{I}$, $x(t)\in{K}$ for every $t\geq{t_0}$, and estimate \cref{eq:21:a,eq:21:b,eq:21:c} holds for every $t\geq{t_0}$. Now we may apply \Cref{Lemma:2} and obtain that estimates \cref{eq:21:d,eq:21:e} hold for every $t\geq{t_0}$. Finally, it follows from estimate \cref{eq:21:a,eq:21:b,eq:21:c}, the choice of $\nu_\eta$, and the variation of constants formula that also estimate \cref{eq:21:f} is satisfied for every $t\geq{t_0}$. This completes the proof of \Cref{Theorem:1}.


\bibliographystyle{abbrv}
\bibliography{bibFile}

\end{document}